\documentclass{asl}

\usepackage[breaklinks=true]{hyperref}
\usepackage{xspace,enumerate}
\usepackage{booktabs}
\usepackage{color,epsfig}
\usepackage{burdges-asl}
\usepackage{calc}

\title{Definably linear groups of finite Morley rank}

\author{Alexandre Borovik}
\address{School of Mathematics\\University of Manchester\\
Oxford Road\\M13 9PL Manchester\\England}

\author{Jeffrey Burdges}
\address{School of Mathematics\\University of Manchester\\
Oxford Road\\M13 9PL Manchester\\England}
\thanks{Supported by NSF postdoctoral fellowship DMS-0503036,
 a Bourse Chateaubriand postdoctoral fellowship, and DFG grant Te 242/3-1.}

\begin{document}

\maketitle

\section*{Introduction}

Zilber's original trichotomy conjecture proposed an explicit classification 
of all Òone-dimensionalÓ objects arising in model theory.  At one point,
classifying the simple groups of finite Morley rank was viewed as a
subproblem whose affirmative answer would justify this conjecture.
Zilber's conjecture was eventually refuted by Hrushovski \cite{Hru_NSMS},
and the classification of simple groups of finite Morley rank
remains open today.  However, these conjectures hold in two significant cases.
First, Hrushovski and Zilber prove the full trichotomy conjecture holds
under very strong geometric assumptions \cite{HZ96}, and this suffices
for various diophantine applications.
Second, the Even \& Mixed Type Theorem \cite{ABC_EvenType} shows that
simple groups of finite Morley rank containing an infinite elementary
abelian 2-subgroup are Chevellay groups over an algebraically closed field
of characteristic two.

In this paper, we clarify some middle ground between these two results
by eliminating involutions from simple groups which are {\em definably}
embedded in a linear group over an algebraically closed field in a structure
of finite Morley rank, and which are not Zariski closed themselves.
One may simplify terminology by saying that $G$ is a {\em definably linear}
group over a field $k$ of finite Morley rank, implicitly using some expansion
of the field language, or just a definably linear group of finite Morley rank.

\begin{namedtheorem}{Theorem}
Let $G$ be a definable infinite simple subgroup of $\GL_n(k)$
 over a field $k$ of finite Morley rank and characteristic zero,
 which is not Zariski closed.
Then $G$ has no involutions, and
 its Borel subgroups are self-normalizing.
\end{namedtheorem}

Poizat has shown that simple groups with such a definable linear embedding
are Chevellay groups when the field has characteristic $p > 0$ \cite{Po01a}.
Much geometric information is available about such groups
 in characteristic zero but not as much as \cite{HZ96} assumes.
For our purposes, the most important geometric fact is that every
element of a counterexample is semisimple in the ambient linear group.

We note that such groups include those groups interpretable in
{\em bad fields}, which have been recently shown to exist \cite{BadField}.

A major feature of our proof is the elimination of
 2-tori of outer automorphisms of simple groups
via the Delehan-Nesin argument
 \cite[\qProposition 13.4]{BN} (see Lemma \ref{bad-like_Sylow_centralize}).
Such an application is a hopeful sign for the current ``$L^*$'' project of
classifying simple groups of finite Morley rank which have an involution.
On the other hand, our proof is driven primarily by one unreasonably
strong fact about these linear groups:
 all strongly real elements, i.e.\ all inverted elements, have
 a nontrivial power inside a unique conjugacy class of Borel subgroups
 which are pairwise disjoint.

In this vein, we leave the proof of the following consequence
 as an exercise to the reader (see also \cite[4.19]{AC03}).

\begin{namedtheorem}{Corollary}
A definably linear group of finite Morley rank is an $L$-groups.
\end{namedtheorem}

In the first section,
 we reduce the problem first to characteristic zero using a result of Poizat,
 and then to the hypotheses used later by our inductive argument.
In the second section,
 we recall various facts about groups without unipotent torsion.
In the third section,
 we prove our main result using the earlier reductions.

The authors thank Oliver \Frecon for helpful comments and corrections.

\section{Definably linear groups}\label{sec:def-linear}

In this section, we strengthen our hypotheses by recalling
 some past work on definably linear groups of finite Morley rank.

Poizat has already shown that simple linear groups over fields
 of characteristic $p > 0$ are isomorphic to algebraic groups.

\begin{fact}[{Poizat \cite[Theorem 1]{Po01a}}]\label{Po01a1}
Let $k$ be a field of finite Morley rank and characteristic $p > 0$.
Then any simple definable subgroup of $\GL_n(k)$ is isomorphic to
 an algebraic group over some field of characteristic $p$,
 but not necessarily $k$.
\end{fact}

As such, our theorem considers only fields of characteristic zero.

\begin{hypothesis}\label{def-linear_hypothesis}
Let $G$ be a definable connected subgroup of $\GL_n(k)$
 over a field $k$ of finite Morley rank and characteristic zero.
\end{hypothesis}

We now show that all definable subgroups satisfy conjugacy of
 Borel subgroups, i.e. maximal connected definable solvable subgroups,
 with the following two facts.

\begin{fact}[{Poizat \cite[Theorem 3]{Po01a}}]\label{Po01a3}
If such a group $G$ is simple, but is not Zariski closed,
 then all elements of $G$ are semisimple.
\end{fact}

\begin{lemma}[{Mustafin \cite[Proposition 2.10]{Mustafin04}}]\label{def-linear_conjugacy}
If our group $G$ consists entirely of semisimple elements,
 then all Borel subgroups of $G$ are abelian and conjugate,
 and this conjugacy class is generic.
\end{lemma}




\begin{proof} 
Let $B$ be a Borel subgroup of $G$.
The Zarisk closure $\hat B$ of $B$ is a connected solvable group.
So its derived subgroup $\hat B'$ consists entirely of unipotent elements.
As $B$ consists entirely of semisimple elements,
 we have $B' \leq B \cap \hat B' = 1$, and $B$ is abelian.
It follows that $\hat B$ is an algebraic torus.

Let $T$ be a maximal algebraic torus of $\GL_n(k)$ which contains $B$.
There are only finitely many distinct intersections
 $T \cap T^g$ for $g\in \GL_n(k)$, by \cite[Rigidity II]{Ch05}.
As any distinct conjugates of $B$ lie in distinct conjugates of $T$,
 there are only finitely many distinct intersections
 $B \cap B^g$ for $g\in G$ too,
and so the union of such intersections is not generic in $B$.
By the genericity argument \cite[Lemma 4.1]{BBC},
 $\bigcup B^G$ is generic in $G$.
Now $G$ has only one conjugacy class of Borel subgroups,
 by \cite{Ja06}. 
\end{proof}

\begin{proposition}\label{def-linear_main}
If $G$ is simple, but not Zariski closed, then
 Borel subgroups of $G$ are abelian and hereditarily conjugate.
\end{proposition}

\begin{proof}
By Fact \ref{Po01a3}, all elements of $G$ are semisimple.
As this property passes to subgroups,
 the result follows from Lemma \ref{def-linear_conjugacy}.
\end{proof}

It follows that our $G$ satisfies Hypothesis \ref{bad-like_hypothesis} below.
So then Theorem \ref{bad-like_main} will show that $G$ has no involutions.

We recall that a Carter subgroup of $G$ is an almost self-normalizing
 definable nilpotent subgroup of $G$
In our case, these are merely our Borel subgroups.
As the Borel subgroups are conjugate,
 the Borel subgroups of $G$ will then be self-normalizing too,
 by Lemma \ref{Weyl_group_cor}.

So these two results will complete the proof.

\section{No unipotent torsion}

We now recall several convenient results from \cite{BC08a}
 concerning groups without {\em unipotent torsion}, i.e.\ which
 contains no definable connected nilpotent subgroup of bounded exponent.

\begin{lemma}\label{pElements_cor}
Let $G$ be a connected group of finite Morley rank without unipotent torsion.
Then every element has a nontrivial power which lies inside a Borel subgroup.
\end{lemma}

\noindent This is a corollary of the following fact.

\begin{fact}[{\cite[Theorem 3]{BC08a}}]\label{pElements}
Let $G$ be a connected group of finite Morley rank of $\pi^\perp$-type,
 i.e. without unipotent $p$-torsion for $p \in \pi$. 
Then all $\pi$-elements of $G$ are toral, i.e.\ lie inside $\pi$-tori.
\end{fact}

We recall that, in this setting,
 divisible torsion groups are referred to as tori and
 their elements may be referred to as toral.
As usual, $H^\o$ denotes the connected component of $H$,
 so a {\em Sylow\o 2-subgroup} is just a maximal 2-torus.

\begin{proof}[Proof of Lemma \ref{pElements_cor}]
Consider an element $a \in G$.
If $d(a)$ is infinite, we choose $n = \abs{d(a)/d(a)^\o}$.
If $d(a)$ is finite, $a$ lies inside some torus by Fact \ref{pElements}.
\end{proof}

We observe that, if one assumes that Borel subgroups are conjugate,
 one may instead prove that all elements lie inside Borel subgroups.
For this, we require that the group has {\em degenerate type}, meaning
 a finite Sylow 2-subgroup or equivalently no involutions \cite[Theorem 1]{BBC}.

\begin{lemma}\label{Weyl_group_cor}
Let $G$ be a connected degenerate type group of finite Morley rank
 without unipotent torsion.
Then $G$ has some self-normalizing Carter subgroup.
\end{lemma}

\noindent This is a corollary of the following.

\begin{definition}
Let $G$ be a group of finite Morley rank, and
 $T$ a maximal divisible abelian torsion subgroup of $G$.
The {\em Weyl group} of $G$ is the finite group $N(T)/C^\o(T)$,
 which can be viewed as a group of automorphisms of $T$.
\end{definition}

\begin{fact}[{\cite[\qTheorem 5]{BC08a}}]\label{Weyl_group}
Let $G$ be a connected group of finite Morley rank.
Suppose the Weyl group is nontrivial and has odd order,
 with $r$ the smallest prime divisor of its order.
Then $G$ contains a unipotent $r$-subgroup.
\end{fact}

\begin{fact}[{\cite[\qLemma 1.2]{BC08a}}]\label{central_maxtorus}
Let $G$ be a connected group of finite Morley rank, 
$T$ a maximal $p$-torus in $G$, and suppose that $T$ is central in $G$.
If $G$ has $p^\perp$-type,
 then any $p$-element $a \in G$ belongs to $T$.
\end{fact}

\begin{proof}[Proof of Lemma \ref{Weyl_group_cor}]
We may assume that $G$ contains torsion.
Let $T$ be a maximal divisible abelian torsion subgroup of $G$.
By \cite{FrJa04},
 $G$ has a Carter subgroup $Q$ containing $T$.
By Fact \ref{Weyl_group}, $N_G(T) = C^\o_G(T)$.
So $N_G(Q) \leq C^\o_G(T)$.
By Fact \ref{central_maxtorus},
 $C^\o_G(T)/d(T)$ is torsion free.
So $N_G(Q) \leq Q$, as desired.
\end{proof}

\section{Bad-like groups}

Here we analyze groups satisfying the following hypotheses.

\begin{hypothesis}\label{bad-like_hypothesis}
Let $G$ be a connected group of finite Morley rank
 without unipotent torsion
whose Borel subgroups are nilpotent and hereditarily conjugate,
 i.e.\ any subgroup $H$ of $G$ satisfies
 $H$-conjugacy of Borel subgroups of $H$.
\end{hypothesis}

These hypotheses clearly pass to subgroups.
We observe that they also pass to sections by nilpotent kernels.

\begin{lemma}\label{bad-like_sections}
Let $K$ be a nilpotent normal subgroup of $G$.
Then $G/K$ satisfies Hypothesis \ref{bad-like_hypothesis} too,
 and images and inverse images preserve Borel subgroups.
\end{lemma}

\begin{proof}
Quotients and extensions of solvable groups are solvable,
so images and inverse images preserve maximal solvable subgroups.
As images also preserve connectivity, they preserve Borel subgroups,
 and Borel subgroups of $G/K$ are nilpotent.
If $K$ is connected, inverse images are connected too,
 and hence preserve Borel subgroups too.
We may now assume that $K$ is finite,
 and even cyclic of order $p^n$.
As $H$ is connected, $K$ is central in $G$.
So, by Fact \ref{pElements},
 $K$ is contained in every maximal $p$-torus of $G$.
By conjugacy, $K$ is contained in every Borel subgroup of $G$ too.
So all inverse images are connected,
 and hence preserve Borel subgroups.
Our first part now follows easily by applying
 the second inside all definable connected subgroups of $G$
\end{proof}

A variation on the Delehan-Nesin argument \cite[\qProposition 13.4]{BN}
 shows that 2-tori have no interesting actions inside such a group.

\begin{proposition}\label{bad-like_Sylow_centralize}
Let $H$ be a definable normal subgroup of $G$ which has no involutions.
Then every 2-torus in $G$ centralizes $H$.
\end{proposition}

For this, we note several facts about such a degenerate type group $H$
 normalized by an involution.

\begin{fact}\label{f:a_nonconj_ainv}
No element $a \in H$ is $H$-conjugate to it's inverse $a^{-1}$.
\end{fact}

\begin{proof}
Suppose that $a^h = a^{-1}$ for some $a,h \in H$.
As $a$ is not an involution,
 we have $h \in N_H(d(a)) \setminus C_H(a)$.
But clearly $h$ has order two modulo $C_H(a)$.
So $H$ contains an involution
 by \cite[Ex.~11 p.~93; Ex.~13c p.~72]{BN}, 
 a contradiction.
\end{proof}

\begin{lemma}\label{l:no_x_in_Ci}
No element of $C_H(i)$ is inverted by any involution in $\gen{i} H$.
\end{lemma}

\begin{proof}
Any involution of $\gen{i} H$ is conjugate to $i$ under $H$.
So otherwise there is an $x \in X^-$ such that $x^h \in C(i)$.
We may assume that $h \in X^-$ too
 because \cite[Ex.\ 14 p.\ 73]{BN} says
 $H = X C_H(i)$ where $X := \{ g\in H \mid g^i = g^{-1} \}$.
So $x^h = (x^h)^i = h x^{-1} h^{-1} = (x^{-1})^{h^{-1}}$, then
 $x^{h^2} = x^{-1}$, contradicting Fact \ref{f:a_nonconj_ainv}.
\end{proof}

\begin{proof}[Proof of Proposition \ref{bad-like_Sylow_centralize}]
Consider a 2-torus $T$ of $G$.  We may take $G := H d(T)$.
We may assume $G$ is centerless, by Lemma \ref{bad-like_sections}. 
So there is a nontrivial strongly real element $i j$ inside $H$.
It follows, by Lemma \ref{pElements_cor}, that some
 nontrivial power $x$ of $i j$ lies inside some Borel subgroup $B$ of $H$.
As all Borel subgroups of $G$ are conjugate, some conjugate $T^h$ of $T$
 lies inside the Borel subgroup of $G$ containing $B$.
But now $T^h$ centralizes $B$ by \cite[Theorem 6.16]{BN}, 
 contradicting Lemma \ref{l:no_x_in_Ci}.
\end{proof}

We say a subgroup $K$ of a group $H$ is {\em subnormal} if
 there is a finite chain
 $$ K \normal K_1 \normal \cdots \normal K_k \normal H\mathperiod $$
A quasisimple subnormal subgroup of a group $H$ is referred to as
 a {\em component} of $H$.
If $H$ is connected, then its components are 
 definable, connected, and normal in $H$, by \cite[Lemma 7.10]{BN}.
They also centralize one another by \cite[Lemma 7.12]{BN}.
We let $E(H)$ denote the central product of the components in $H$,
which happen to be connected, and set $F^*(H) := F^\o(H) E(H)$.
According to the following,
 $\Aut(F^*(\bar{G}))$ controls the structure of $\bar{G}$.

\begin{fact}[{\cite[\qTheorem 7.13 \& \qCorollary 7.14]{BN}}]\label{genfitting}
Any group $H$ of finite Morley rank satisfies
$$ C_H(F^*(H)) = Z(F^*(H)) \leq F^*(H)
    \quad\text{and}\quad
    H/Z(F^*(H)) \leq \Aut(F^*(H))\mathperiod $$
\end{fact}

\begin{lemma}\label{bad-like_Dgrp}
If all simple sections of $G$ with nilpotent kernels have degenerate type,
 then $G$ has a unique Sylow 2-subgroup, which is connected and central.
\end{lemma}

\begin{proof}
A connected degenerate type group of finite Morley rank
 has no involutions by \cite[Theorem 1]{BBC}.
So, as $E(G)/Z(E(G))$ is $2^\perp$ by hypothesis,
 we find that $E(G)$ is $2^\perp$ too.
By Lemma \ref{bad-like_Sylow_centralize},
 the Sylow\o 2-subgroup $S^\o$ of $G$ centralizes $E(G)$.
But $G \hookrightarrow \Aut(E(G))$ by Fact \ref{genfitting}.
So any 2-torsion in $G$ is lies inside the solvable radical $\sigma(G)$,
 which lies inside $F(G)$ because Borel subgroups are nilpotent.
It follows that the Sylow 2-subgroup is connected
 by Fact \ref{pElements} and \cite[Theorem 9.29]{BN},
 and central by \cite[Theorem 6.16]{BN}.
\end{proof}

We now prove our main result.

\begin{theorem}\label{bad-like_main}
Suppose that $G$ is simple.
Then $G$ has no involutions.
\end{theorem}

\begin{proof}
We consider a counterexample $G$ of minimal Morley rank.
So, by Lemma \ref{bad-like_sections},
 all proper simple section with nilpotent kernels have degenerate type.
By Lemma \ref{bad-like_Dgrp},
 any proper connected subgroup $H$ centralizes its Sylow 2-subgroup,
 which is a 2-torus.
As $G$ is simple, it follows that 
\begin{texteqn}{\tag{$\star_1$}}
the centralizers of distinct Sylow\o 2-subgroups are disjoint.
\end{texteqn}
Also, any Sylow\o 2-subgroup is central in a Borel subgroup.
By conjugacy, all Borel subgroups have a central Sylow\o 2-subgroup.
So
\begin{texteqn}{\tag{$\star_2$}}
 two Borel subgroups meet iff they both
 contain the same Sylow\o 2-subgroup of $G$.
\end{texteqn}

\begin{Claim}\label{bad-like_involution_pairs}
Any two distinct involutions $i$ and $j$ normalize
 a common Sylow\o 2-subgroup $T$ of $G$, and $C^\o_G(T)$.
\end{Claim}

\begin{verification}
By Lemma \ref{pElements_cor},
 there is a nontrivial power $x$ of $i j$ inside a Borel subgroup $Q$ of $G$.
As $i$ inverts $x$, $Q^i$ also contains $x$.
By $\star_2$,
 $Q$ and $Q^i$ contain the same Sylow\o 2-subgroup $T$ of $G$.
So $i$ normalizes $T$.
Similarly, $j$ normalizes $T$ too.
\end{verification}

Every involution is contained in a Sylow\o 2-subgroup,
 by Fact \ref{pElements}, which is unique by $\star_1$.
Let $T_i$ denote the Sylow\o 2-subgroup containing $i \in I(G)$.
Now $C^\o_G(i) = C^\o_G(T_i)$ since
 $T_i$ is the unique Sylow\o 2-subgroup of $C^\o_G(i)$, by $\star_1$.

\begin{Claim}\label{bad-like_inverted_Borel}
Let $T$ be a Sylow\o 2-subgroup of $G$,
 and let $j$ be an involution normalizing $T$.
Then either $j\in T$ or $j$ inverts $C_G(T)$.
In particular, $C^\o_G(T)$ is an abelian Borel subgroup of $G$.
\end{Claim}

\begin{verification}
Suppose that $j \notin T$.
Then $C_G(T,j) = C_G(T_j) \cap C_G(T) = 1$.
So, by \cite[Ex.\ 13 \& 15, p.\ 79]{BN},
 $C^\o_G(T)$ is abelian and inverted by $j$.
Now $C^\o_G(T)$ is a Borel subgroup of $G$,
 as it containes a Borel subgroup of $G$.
\end{verification}

As a consequence, if involutions $i$ and $j$ do not commute,
 then neither $i$ nor $j$ lie inside a Sylow\o 2-subgroup which
 they both normalize, and both invert its Borel subgroup.

\begin{Claim}\label{bad-like_pr_1}
$G$ has \Prufer rank one.
So a Sylow 2-subgroup of $G$ is isomorphic to that of $\PSL_2$.
\end{Claim}

\begin{verification}
Suppose towards a contradiction that $\pr(G) > 1$.
Choose a pair of non-commuting involutions $i$ and $j$,
 and fix a Sylow\o 2-subgroup $T$ which is normalized by them both.
As $i$ inverts $C^\o_G(T)$ by Claim \ref{bad-like_inverted_Borel},
 $\Omega_1(T)$ normalizes $C^\o_G(T_i)$.
By generation \cite[Theorem 1]{BC08a},
 $C^\o_G(T_i) = \gen{ C^\o_G(T_i,u) \mid u\in I(T) }$.
So $C^\o_G(T_i,u) \neq 1$ for some $u\in I(T)$,
 which yields the contradiction $T = T_i$.
Thus $\pr(G) = 1$.

Let $S$ be a Sylow 2-subgroup of $G$.
By Fact \ref{pElements},
 $C_S(S^\o) = S^\o$, and $[S,S^\o] \leq 2$.
There are involutions in $S \setminus S^\o$
 by Claims \ref{bad-like_involution_pairs} and \ref{bad-like_inverted_Borel}.
So $S$ has the desired structure.
\end{verification}

Consider the set $I$ of all involutions in $G$.
We define the set $L$ of lines on $I$ to be the collection of
cosets $j Q$ where $Q$ is a Borel subgroup of $G$ and $j$ inverts $Q$.
We now prove that $I$ and $L$ form a projective plane.


\begin{namedtheorem}{PP2}
Any two points lie on a unique line.
\end{namedtheorem}

Any two involutions $i,j \in I$ normalize some Sylow\o 2-subgroup $T$,
 by Claim \ref{bad-like_involution_pairs}.
By Claim \ref{bad-like_inverted_Borel},
 $Q := C^\o_G(T)$ is a Borel subgroup inverted by $i$ and $j$.
So $i$ and $j$ lie on the line $j Q$.
Let $x$ denote the power of $i j$ lies inside some $Q^g$.
By $\star_1$, $T^g$ is the unique Sylow\o 2-subgroup of $C^\o_G(x)$.
So the line $j Q$ is unique.

\begin{namedtheorem}{PP3}
Any two lines intersect at a unique point.
\end{namedtheorem}

For any two lines $i_1 Q_1$ and $i_2 Q_2$,
 we consider the involutions $j_1 \in I(Q_1)$ and $j_2 \in I(Q_2)$,
 which are unique by Claim \ref{bad-like_pr_1}.
Again, there is a Borel subgroup $R$ inverted by both,
 by Claims \ref{bad-like_involution_pairs} and \ref{bad-like_inverted_Borel}.
As $j_1$ and $j_2$ both centralize the involution $t \in I(R)$,
 $t$ inverts both $Q_1$ and $Q_2$,
 by Claims \ref{bad-like_inverted_Borel} and \ref{bad-like_pr_1}.
Thus the lines $i_1 Q_1$ and $i_2 Q_2$ meet at $t$.
Such a point is unique by PP2.

\begin{namedtheorem}{PP4}
There are four points, no three of which are colinear.
\end{namedtheorem}

Consider a Borel subgroup $Q$ of $G$, and
 the unique involution $i \in I(Q)$. 
Let $x$, $y$, and $z$ be three involutions inverting $Q$,
 and let $j$ be a third point on the line through $i$ and $z$.
The original points $i,x,y$ are obviously not colinear.
The pair $i,j$ is not colinear with either $x$ or $y$ by PP3.
Also $j,x,y$ is not colinear by PP2.
So $i,j,x,y$ are the four desired involutions.

We prove one additional property of this projective plane.

\begin{namedtheorem}{PP+}
Any three involutions $x,y,z\in I(G)$ are colinear iff
 their product $x y z$ is an involution.
\end{namedtheorem}

The product of three colinear involutions is an involution
 by Claim \ref{bad-like_pr_1}.
Consider three involutions $x,y,z\in I(G)$ whose
 product $x y z$ is an involution $i$.
Let $Q$ be the Borel subgroup such that $y$ and $z$ invert $Q$.
As $x y z$ is an involution, $(y z)^x = (y z)^{-1}$. 
So $x$ normalizes $Q$ by $\star_1$.
But $x \notin Q$ as $Q$ has only one involution. 
So $x$ inverts $Q$, and $x$ lies on the line with $y$ and $z$.

\smallskip

These four conditions contradict the following theorem. 

\begin{namedtheorem}{Bachmann's Theorem}[{\cite[Theorems 8.15 \& 8.18]{BN}}]
Let $G$ be a group whose set of involutions $I$ posses the structure of
a projective plane, and three involutions are colinear iff their product
is an involution.  Then $\gen{I} \cong \SO_3(k,f)$ for some
interpretable field $k$ and some non-isotropic quadratic form $f$.
In particular, $G$ does not have finite Morley rank.
\end{namedtheorem}

This concludes the proof of Theorem \ref{bad-like_main}.
\end{proof}

\smallskip

We conclude by conjecturing a a stronger variation the above style of argument
 exploiting Bachmann's Theorem.

\begin{conjecture}
There is no simple group of finite Morley rank in which
 every strongly real elements lies inside some $C^\o(T)$
 for $T$ is a Sylow\o 2-subgroup.
\end{conjecture}


\small
\bibliographystyle{asl}
\bibliography{burdges,fMr}

\end{document}